\documentclass[a4paper,10pt]{article}
\usepackage{stmaryrd}
\usepackage{amsfonts}
\usepackage{bbm}
\usepackage{amscd}
\usepackage{mathrsfs}
\usepackage{latexsym,amssymb,amsmath,amscd,amscd,amsthm,amsxtra}
\usepackage[dvips]{graphicx}
\usepackage[utf8]{inputenc}
\usepackage[T1]{fontenc}
\usepackage{lmodern}
\usepackage{amssymb}
\usepackage[all]{xy}
\usepackage{nicefrac,mathtools,enumitem}
\usepackage{microtype}
\usepackage{xcolor}
\newtheorem{thm}{Theorem}[section]
\newtheorem{lem}[thm]{Lemma}
\newtheorem{cor}[thm]{Corollary}
\newtheorem{pro}[thm]{Proposition}

\newtheorem{rmk}[thm]{Remark}
\newtheorem{defi}[thm]{Definition}

\newcommand {\emptycomment}[1]{}

\newcommand{\lon }{\,\rightarrow\,}
\newcommand{\be }{\begin{equation}}
\newcommand{\ee }{\end{equation}}
\newcommand{\add}{\frka\frkd}

\newcommand{\pf}{\noindent{\bf Proof.}\ }

\newcommand{\g}{\frkg}
\newcommand{\h}{\frkh}

\newcommand{\huaB}{\mathcal{B}}

\newcommand{\huaH}{\mathcal{H}}

\newcommand{\huaO}{\mathcal{O}}

\newcommand{\frka}{\mathfrak a}

\newcommand{\frkd}{\mathfrak d}

\newcommand{\frkg}{\mathfrak g}
\newcommand{\frkh}{\mathfrak h}

\def\qed{\hfill ~\vrule height6pt width6pt depth0pt}

\newcommand{\br}[1]{   [ \cdot,    \cdot  ]_\frkg   }

\newcommand{\Id}{\rm{Id}}

\newcommand{\HF}{\mathsf{HF}}
\newcommand{\Der}{\mathsf{Der}}

\newcommand{\Inn}{\mathsf{Inn}}
\newcommand{\Out}{\mathsf{Out}}
\newcommand{\Ad}{\mathsf{Ad}}

\newcommand{\gl}{\mathfrak {gl}}

\newcommand{\ad}{\mathsf{ad}}
\newcommand{\pr}{\mathrm{pr}}

\newcommand{\K}{\mathbb{K}}

\begin{document}
\title{
{ Cohomologies, deformations  and extensions of $n$-Hom-Lie algebras}
\thanks
 {
Research supported by NSFC (11471139) and NSF of Jilin Province (20170101050JC).
 }
}
\author{Lina Song and Rong Tang  \\
Department of Mathematics, Jilin University,\\
 Changchun 130012, Jilin, China
\\\vspace{3mm}
Email: songln@jlu.edu.cn, tangrong16@mails.jlu.edu.cn}

\date{}


\maketitle
\begin{abstract}
In this paper, first we give the cohomologies of an $n$-Hom-Lie algebra and  introduce the notion of a derivation of an $n$-Hom-Lie algebra.    We show that a derivation of an $n$-Hom-Lie algebra is a $1$-cocycle with the coefficient in the adjoint representation. We also give the formula of the dual representation of a representation of an $n$-Hom-Lie algebra. Then, we study $(n-1)$-order deformation of an $n$-Hom-Lie algebra. We introduce the notion of a Hom-Nijenhuis operator, which could generate a trivial  $(n-1)$-order deformation of an $n$-Hom-Lie algebra. Finally, we introduce the notion of a generalized derivation of an $n$-Hom-Lie algebra, by which we can construct a new $n$-Hom-Lie algebra, which is called the generalized derivation extension of an $n$-Hom-Lie algebra.
\end{abstract}

\tableofcontents

\section{Introduction}

The representation theory of an algebra is very important to this algebraic structure. Given a representation, one can obtain the corresponding cohomology, which could provide invariants. The cohomology plays important roles in the study of deformation and extension problems.

The notion of a Hom-Lie algebra was introduced by Hartwig, Larsson, and Silvestrov in \cite{HLS} as part of a study of deformations of the Witt and the Virasoro algebras. In a Hom-Lie algebra, the Jacobi identity is twisted by a linear map, called the Hom-Jacobi identity. Some $q$-deformations of the Witt and the Virasoro algebras have the structure of a Hom-Lie algebra \cite{HLS,hu}. Because of close relation to discrete and deformed vector fields and differential calculus \cite{HLS,LD1,LD2}, more people pay special attention to this algebraic structure. In particular, representations and cohomologies of Hom-Lie algebras are studied in \cite{AEM,MS1,sheng}.
On the other hand, the notion of an $n$-Hom-Lie algebra was introduced in \cite{Hom-n-Lie}, which is a generalization of an $n$-Lie algebra introduced in \cite{Filippov}. See the review article \cite{review} for more information about $n$-Lie algebras. Then several aspects about $n$-Hom-Lie algebras are studied. For example, the cohomologies adapted to central extensions and deformations are studied in \cite{Ammarb-Mabroukb-Makhlouf}; 2-cocycles that used to studied abelian extensions are studied in \cite{bai}; Construction of 3-Hom-Lie algebras from Hom-Lie algebras are studied in \cite{AMS}, and extensions of 3-Hom-Lie algebras are studied in \cite{chen}.  However, the systematic study of the cohomology of an $n$-Hom-Lie algebra is still lost.

The first purpose of this paper is to give a systematic study of cohomology of $n$-Hom-Lie algebras. We introduce the coboundary operator associated to a general representation of an $n$-Hom-Lie algebra and obtain the corresponding cohomology. As a byproduct, we introduce the notion of a derivation of an $n$-Hom-Lie algebra, which is different from existing ones. We show that a derivation is exactly a 1-cocycle with the coefficient in the adjoint representation, which generalizes the fact about derivations of a Lie algebra. We also studied $(n-1)$-order deformations of $n$-Hom-Lie algebras, which is inspired by the work in \cite{LSZB}. It turns out that for $n$-ary algebras, one should study $(n-1)$-order deformations, instead of 1-order deformations, to obtain invariants. The second purpose is to give the correct definition of a dual representation of a representation of an $n$-Hom-Lie algebra. In \cite{chen}, the authors claimed that  the usual definition of $\ad^*$ is still a representation without a proof. To solve this problem, in \cite{bai} the authors add a strong condition on a representation $\rho$ of an $n$-Hom-Lie algebra to make $\rho^*$ to be a representation, which generalizes the idea from \cite{sheng1}. Here we give a new formula to define the dual representation, which is quite natural. The last purpose is to give an approach to construct new $n$-Hom-Lie algebras. To do this, we introduce the notion of a generalized derivation of an $n$-Hom-Lie algebra, by which we can construct a new $n$-Hom-Lie algebra, which is called generalized derivation extension.

The paper is organized as follows. In Section 2, we recall some basic notions of Hom-Leibniz algebras and representations of $n$-Hom-Lie algebras. In Section 3, we study cohomologies, derivations and dual representations of $n$-Hom-Lie algebras. In Section 4, we studied $(n-1)$-order deformations of $n$-Hom-Lie algebras and introduce the notion of Hom-Nijenhuis operators which could generate trivial deformations. In Section 5, we introduce the notion of a generalized derivation of an $n$-Hom-Lie algebra, and construct a new $n$-Hom-Lie algebra, which is called generalized derivation extension.

In this paper, we work over an algebraically closed field $\K$ of characteristic 0 and all the vector spaces are over $\K$.

\vspace{2mm}
 \noindent {\bf Acknowledgement:} We give our warmest thanks to Yunhe Sheng  for very helpful suggestions that improve the paper.

\section{Preliminaries}

In this section, we recall Hom-Leibniz algebras and representations of $n$-Hom-Lie algebras.

\begin{defi}
  A (multiplicative) {\bf Hom-Leibniz algebra} is a vector space $\g$ together with a bracket operation $[\cdot, \cdot]_\g:\g\otimes\g\longrightarrow\g$ and an algebraic automorphism   $\alpha:\g\longrightarrow\g$, such that for all $x,y,z\in\g$, we have
  \begin{equation}
    [\alpha(x),[y,z]_\g]_\g=[[x,y]_\g,\alpha(z)]_\g+[\alpha(y),[x,z]_\g]_\g.
  \end{equation}
\end{defi}
In particular, if the bracket operation $[\cdot,\cdot]_\g$ is skew-symmetric, we obtain the definition of a Hom-Lie algebra.

Let $V$ be a vector space, and $\beta\in GL(V)$. Define a skew-symmetric bilinear bracket operation $[\cdot,\cdot]_{\beta}:\wedge^2\mathfrak{gl}(V)\longrightarrow\mathfrak{gl}(V)$ by
\begin{eqnarray}\label{eq:bracket}
[A,B]_{\beta}=\beta \circ A \circ\beta^{-1}\circ B \circ\beta^{-1}-\beta\circ B \circ\beta^{-1}\circ A\circ \beta^{-1}, \hspace{3mm}\forall A,B\in \mathfrak{gl}(V).
\end{eqnarray}
Denote by $\Ad_{\beta}:\mathfrak{gl}(V)\lon \mathfrak{gl}(V)$ the adjoint action on $\mathfrak{gl}(V)$, i.e.
\begin{equation}\label{eq:Ad}
\Ad_{\beta}(A)=\beta\circ A\circ \beta^{-1}.
\end{equation}
\begin{pro}{\rm (\cite[Proposition 4.1]{shengxiong})}\label{pro:Hom-Lie}
   With the above notations, $(\mathfrak{gl}(V),[\cdot,\cdot]_{\beta},\Ad_{\beta})$ is a regular Hom-Lie algebra.
 \end{pro}This Hom-Lie algebra plays an important role in the representation theory of Hom-Lie algebras. See \cite{shengxiong} for more details.

\begin{defi}
   An $n$-Hom-Lie algebra is a vector space $\g$ equipped with a bracket operation $[\cdot,\cdots,\cdot]_\g:\wedge^n\g\longrightarrow\g$ and an algebraic automorphism   $\alpha:\g\longrightarrow\g$ such that for all $x_1,\cdots,x_{n-1}$, $y_1,\cdots, y_n\in\g$, the following {\bf Hom-Fundamental} identity  holds:
\begin{eqnarray}
\nonumber \HF_{x_1,\cdots,x_{n-1},y_1,\cdots,y_n}
&\triangleq&[\alpha(x_1),\cdots,\alpha(x_{n-1}),[y_1,\cdots,y_n]_\g]_\g\\\nonumber&&-\sum_{i=1}^n[\alpha(y_1),\cdots,\alpha(y_{i-1}),
     [x_1,\cdots,x_{n-1},y_i]_\g,\alpha(y_{i+1}),\cdots,\alpha(y_n)]_\g\\
\label{eq:de1}&=&0.
\end{eqnarray}
 \end{defi}

Any linear map $\alpha:\g\longrightarrow\g$ induces a linear map $\widetilde{\alpha}:\wedge
^{n-1}\g\longrightarrow \wedge^{n-1}\g$ via
\begin{equation}
  \widetilde{\alpha}(x_1\wedge\cdots\wedge x_{n-1})=\alpha(x_1)\wedge\cdots\wedge\alpha(x_{n-1}).
\end{equation}
Similar as the case of $n$-Lie algebras, elements in $\wedge^{n-1}\g$ are called fundamental elements.  On $\wedge^{n-1}\g$, one can define a new bracket operation $[\cdot,\cdot]_F$ by
\begin{equation}
  [X,Y]_F=\sum_{i=1}^{n-1}\alpha(y_1)\wedge\cdots\wedge\alpha(y_{i-1})\wedge[x_1,\cdots,x_{n-1},y_i]_\g\wedge\alpha(y_{i+1})
  \wedge\cdots\wedge\alpha(y_{n-1}),
\end{equation}
for all $X=x_1\wedge\cdots \wedge x_{n-1}$ and $Y=y_1\wedge\cdots\wedge y_{n-1}$.
It is proved in \cite{Ammarb-Mabroukb-Makhlouf} that $(\wedge^{n-1}\g,[\cdot,\cdot]_F,\widetilde{\alpha})$ is a Hom-Leibniz algebra.

\begin{defi}
  A {\bf morphism} of $n$-Hom-Lie algebras $f:(\mathfrak{g},[\cdot,\cdots,\cdot]_{\mathfrak{g}},\alpha)\lon (\h,[\cdot,\cdots,\cdot]_\h,\gamma)$ is a linear map $f:\mathfrak{g}\lon \mathfrak{h}$ such that
  \begin{eqnarray}
f[x_1,\cdots,x_n]_{\mathfrak{g}}&=&[f(x_1),\cdots,f(x_n)]_{\mathfrak{h}},\hspace{3mm}\forall x_1,\cdots,x_n\in \mathfrak{g},\\
    f\circ \alpha&=&\gamma\circ f.
  \end{eqnarray}
\end{defi}

\begin{defi}\label{defi:rep}
  A {\bf representation} of an $n$-Hom-Lie algebra $(\g,[\cdot,\cdots,\cdot]_\g,\alpha)$ on a vector space $V$ with respect to a linear automorphism $\beta\in GL(V)$ is a linear map $\rho:\g\longrightarrow\gl(V)$ such that for all $X,Y\in\wedge ^{n-1}\g$ and $x_1,\cdots,x_{n-2},y_1,\cdots,y_n\in\g$, we have
\begin{itemize}
  \item[{\rm(i)}] $\rho(\widetilde{\alpha}(X) )\circ\beta=\beta\circ \rho(X);$

  \item[{\rm(ii)}] $ \rho(\widetilde{\alpha}(X))\circ\rho(Y)-\rho(\widetilde{\alpha}(Y))\circ \rho(X)=\rho([X,Y]_F)\circ \beta;$

    \item[{\rm(iii)}] \begin{eqnarray*}
        &&\rho(\alpha(x_1),\cdots,\alpha(x_{n-2}),[y_1,\cdots,y_n]_\g)\circ\beta\\
        &=&\sum_{i=1}^n(-1)^{n-i}\rho(\alpha(y_1),\cdots,\hat{\alpha}(y_i),\cdots,\alpha(y_n))\circ\rho(x_1,\cdots,x_{n-2},y_i).
        \end{eqnarray*}
\end{itemize}
\end{defi}
We denote a representation by $(V,\rho,\beta)$.
\begin{rmk}
  Condition (i) and (ii) in the above definition means that $\rho:\wedge^{n-1}\g\longrightarrow \gl(\g)$
is a morphism from the Hom-Leibniz algebra $(\wedge^{n-1}\g,[\cdot,\cdot]_F,\widetilde{\alpha})$ to the Hom-Lie algebra $(\mathfrak{gl}(V),[\cdot,\cdot]_{\beta},\Ad_{\beta})$.
\end{rmk}

Define $\ad:\wedge^{n-1}\g\longrightarrow\gl(\g)$ by
\begin{equation}
  \ad_Xy=[x_1,\cdots,x_{n-1},y],\quad\forall X=x_1\wedge\cdots \wedge x_{n-1}\in\wedge^{n-1}\g,~y\in\g.
\end{equation}
Then $(\g,\ad,\alpha)$ is a representation of the $n$-Hom-Lie algebra $(\g,[\cdot,\cdots,\cdot]_\g,\alpha)$ on $\g$ with respect to $\alpha$, which is called the {\bf adjoint representation}.

The following result is straightforward.

\begin{pro}\label{tri}
Let  $(V,\rho,\beta)$ be a representation of an $n$-Hom-Lie algebra $(\g,[\cdot,\cdots,\cdot]_\g,\alpha)$. Define a bracket operation $[\cdot,\cdots,\cdot]_{\rho}:\wedge^n(\g\oplus V)\lon\g\oplus V$ by
$$[x_1+u_1,\cdots,x_n+u_n]_{\rho}=[x_1,\cdots,x_n]_\g+\sum_{i=1}^{n}(-1)^{n-i}\rho(x_1,\cdots,\hat{x}_i,\cdots,x_n)u_i,\,\,\,\,\forall x_i\in\g,u_i\in V.$$
Define $\alpha+\beta:\g\oplus V\lon\g\oplus V$ by
$$(\alpha+\beta)(x+u)=\alpha(x)+\beta(u).$$
Then $(\g\oplus V,[\cdot,\cdots,\cdot]_{\rho},\alpha+\beta)$ is an $n$-Hom-Lie algebra, which we call the {\bf semi-direct product} of the  $n$-Hom-Lie algebra $(\g,[\cdot,\cdots,\cdot]_\g,\alpha)$ by the representation $(V,\rho,\beta)$. We denote this semidirect $n$-Hom-Lie algebra simply by $\g\ltimes V$.
\end{pro}

 \section{Cohomologies, derivations and dual representations of $n$-Hom-Lie algebras}

  \subsection{Cohomologies of $n$-Hom-Lie algebras}
  Let $(V,\rho,\beta)$  be a representation of an $n$-Hom-Lie algebra $(\g,[\cdot,\cdots,\cdot]_\g,\alpha)$.
A $p$-cochain on  $\frkg$ with the coefficients in a representation   $(V,\rho,\beta)$ is  a multilinear  map $$f:\wedge^{n-1}\frkg\otimes\stackrel{(p-1)}{\ldots}\otimes\wedge^{n-1}\frkg\wedge\frkg\longrightarrow V.$$
Denote the space of $p$-cochains by $C^{p}(\g;V).$ Define the coboundary operator $\delta:C^{p}(\g;V)\longrightarrow C^{p+1}(\g;V)$   by
\begin{eqnarray*}
&&(\delta f)(X_1,\ldots,X_p,z)\\
&=& \sum_{1\leq i<k\leq p}(-1)^i\beta f\Big(\widetilde{\alpha}^{-1}(X_1),\cdots,\hat{X}_i,\cdots,\widetilde{\alpha}^{-1}(X_{k-1}),[\widetilde{\alpha}^{-2}(X_i), \widetilde{\alpha}^{-2}(X_k)]_F,\\
&&\widetilde{\alpha}^{-1}(X_{k+1}),\cdots,\widetilde{\alpha}^{-1}(X_{p}),\alpha^{-1}(z)\Big)\\
&&+\sum_{i=1}^p(-1)^i\beta f(\widetilde{\alpha}^{-1}(X_1),\cdots,\hat{X}_i,\cdots,\widetilde{\alpha}^{-1}(X_{p}),
[\widetilde{\alpha}^{-2}(X_i),\alpha^{-2}(z)]_\g)\\
&&+\sum_{i=1}^p(-1)^{i+1}\rho(X_i) f(\widetilde{\alpha}^{-1}(X_1),\cdots,\hat{X}_i,\cdots,\widetilde{\alpha}^{-1}(X_{p}),\alpha^{-1}(z))\\
&&+\sum_{i=1}^{n-1}(-1)^{n+p-i+1}\rho(x^1_p,x^2_p,\cdots,\hat{x^i_p},\cdots,x^{n-1}_p,z)
f(\widetilde{\alpha}^{-1}(X_1),\cdots,\widetilde{\alpha}^{-1}(X_{p-1}),\alpha^{-1}(x^i_p)) ,
\end{eqnarray*}
for all   $X_i=(x^1_i,x^2_i,\cdots,x^{n-1}_i)\in\wedge^{n-1}\frkg$ and $z\in\frkg.$

\begin{thm}
  With the above notations, $\delta^2=0$. Thus, we have a well-defined cohomological complex $\big(C^\bullet(\g;V)=\oplus_{p\ge1} C^{p}(\g;V),\delta\big)$.
\end{thm}
\pf It follows by straightforward computations. We omit details. \qed\vspace{3mm}

A $p$-cochain $f\in C^{p}(\g;V)$ is called a {\bf $p$-cocycle} if $\delta(f)=0$. A $p$-cochain $f\in C^{p}(\g;V)$ is called a {\bf $p$-coboundary} if $f=\delta(g)$ for some $g\in C^{p-1}(\g;V)$. Denote by $\mathcal{Z}^p(\g;V)$ and $\mathcal{B}^p(\g;V)$ the sets of $p$-cocycles and  $p$-coboundaries respectively. We define the $p$-th cohomology group
$\mathcal{H}^p(\g;V)$ to be $\mathcal{Z}^p(\g;V)/\mathcal{B}^p(\g;V)$.

  \subsection{Derivations of $n$-Hom-Lie algebras}

Generalizing the notion of a derivation of a Hom-Lie algebra given in \cite{songtang}, we give the notion of a derivation of an $n$-Hom-Lie algebra as follows.

  \begin{defi}
    A {\bf derivation} of an $n$-Hom-Lie algebra $(\g,[\cdot,\cdots,\cdot]_\g,\alpha)$ is a linear map $D\in\gl(\g)$, such that for all $x_1,\cdots,x_n\in\g$, the following equality holds:
        \begin{equation}\label{defi:derivation}
      D[x_1,\cdots,x_n]_\g=\sum_{i=1}^n[\alpha(x_1),\cdots,\alpha(x_{i-1}),(\Ad_{\alpha^{-1}}D)(x_i),\alpha(x_{i+1}),\cdots,\alpha(x_n)]_\g.
    \end{equation}
  \end{defi}

  Denote the set of derivations by $\Der(\g)$.

  \begin{lem}\label{lem:sub1}
 For all $D\in \Der(\g)$, we have $\Ad_{\alpha}D\in\Der(\g)$.
\end{lem}
\pf For all $x_1,\cdots,x_{n}\in\g$, we have
\begin{eqnarray*}
 \Ad_{\alpha}D[x_1,\cdots,x_n]_\g&=& \alpha D[\alpha^{-1}(x_1),\cdots,\alpha^{-1}(x_n)]_\g
 \\&=&\alpha \Big(\sum_{i=1}^n[x_1,\cdots,x_{i-1},(\Ad_{\alpha^{-1}}D)(\alpha^{-1}(x_i)), x_{i+1},\cdots,x_n]_\g\Big)\\
 &=&\sum_{i=1}^n[\alpha(x_1),\cdots,\alpha(x_{i-1}),(\Ad_{\alpha^{-1}}\Ad_{\alpha}D)(x_i),\alpha(x_{i+1}),\cdots,\alpha(x_n)]_\g,
\end{eqnarray*}
which implies that $\Ad_{\alpha}D$ is also a derivation. \qed

\begin{lem}\label{lem:sub2}
For all $D,D'\in \Der(\mathfrak{g})$, we have
$[D,D']_{\alpha}\in \Der(\mathfrak{g}).$
\end{lem}
\pf For all $x_1,\cdots,x_n\in \g$, by \eqref{eq:bracket} and \eqref{defi:derivation} we have
\begin{eqnarray*}
&&[D,D']_{\alpha}([x_1,\cdots,x_n]_\g)\\&=&(\alpha\circ D\circ \alpha^{-1}\circ D'\circ\alpha^{-1}-\alpha\circ D'\circ \alpha^{-1}\circ D\circ\alpha^{-1})[x_1,\cdots,x_n]_\g\\
&=&(\alpha\circ D\circ\alpha^{-1})\sum_{i=1}^{n}[x_1,\cdots,x_{i-1},(\Ad_{\alpha^{-1}}D')\alpha^{-1}x_i,x_{i+1},\cdots,x_n]_\g\\
&&-(\alpha\circ D'\circ \alpha^{-1})\sum_{i=1}^{n}[x_1,\cdots,x_{i-1},(\Ad_{\alpha^{-1}}D)\alpha^{-1}x_i,x_{i+1},\cdots,x_n]_\g\\
&=&(\alpha\circ D)\sum_{i=1}^{n}[\alpha^{-1}x_1,\cdots,\alpha^{-1}x_{i-1},\alpha^{-2}D'x_i,\alpha^{-1}x_{i+1},\cdots,\alpha^{-1}x_n]_\g\\
&&-(\alpha\circ D')\sum_{i=1}^{n}[\alpha^{-1}x_1,\cdots,\alpha^{-1}x_{i-1},\alpha^{-2}Dx_i,\alpha^{-1}x_{i+1},\cdots,\alpha^{-1}x_n]_\g\\
&=&\alpha\sum_{i<j}[x_1,\cdots,x_{i-1},\alpha^{-1}D'x_i,x_{i+1},\cdots,x_{j-1},(\Ad_{\alpha^{-1}}D)\alpha^{-1}x_j,x_{j+1},\cdots,x_n]_\g\\
&&+\alpha\sum_{i=1}^{n}[x_1,\cdots,x_{i-1},(\Ad_{\alpha^{-1}}D)\alpha^{-2}D'x_i,x_{i+1},\cdots,x_n]_\g\\
&&+\alpha\sum_{i>j}[x_1,\cdots,x_{j-1},(\Ad_{\alpha^{-1}}D)\alpha^{-1}x_j,x_{j+1},\cdots,x_{i-1},\alpha^{-1}D'x_i,x_{i+1},\cdots,x_n]_\g\\
&&-\alpha\sum_{i<j}[x_1,\cdots,x_{i-1},\alpha^{-1}Dx_i,x_{i+1},\cdots,x_{j-1},(\Ad_{\alpha^{-1}}D')\alpha^{-1}x_j,x_{j+1},\cdots,x_n]_\g\\
&&-\alpha\sum_{i=1}^{n}[x_1,\cdots,x_{i-1},(\Ad_{\alpha^{-1}}D')\alpha^{-2}Dx_i,x_{i+1},\cdots,x_n]_\g\\
&&-\alpha\sum_{i>j}[x_1,\cdots,x_{j-1},(\Ad_{\alpha^{-1}}D')\alpha^{-1}x_j,x_{j+1},\cdots,x_{i-1},\alpha^{-1}Dx_i,x_{i+1},\cdots,x_n]_\g\\
&=&\sum_{i=1}^{n}[\alpha(x_1),\cdots,\alpha(x_{i-1}),(\Ad_{\alpha^{-1}}[D,D']_{\alpha})(x_i),\alpha(x_{i+1}),\cdots,\alpha(x_n)]_\g.
\end{eqnarray*}
Therefore, we have $[D,D']_{\alpha}\in \Der(\mathfrak{g})$.\qed

\begin{pro} \label{Der}
 With the above notations, $(\Der(\g),[\cdot,\cdot]_{\alpha},\Ad_{\alpha})$ is a Hom-Lie algebra, which is a subalgebra of the Hom-Lie algebra $(\mathfrak{gl}(\frkg),[\cdot,\cdot]_{\alpha},\Ad_{\alpha})$ given in Proposition \ref{pro:Hom-Lie}.
 \end{pro}
 \pf By Lemma \ref{lem:sub1} and \ref{lem:sub2}, $(\Der(\g),[\cdot,\cdot]_{\alpha},\Ad_{\alpha})$ is a Hom-Lie subalgebra of the Hom-Lie algebra $(\mathfrak{gl}(\frkg),[\cdot,\cdot]_{\alpha},\Ad_{\alpha})$.  \qed

  \begin{pro}
Let $(\g,[\cdot,\cdots,\cdot]_\g,\alpha)$ be an $n$-Hom-Lie algebra. $f\in C^1(\g;\g)$ is a $1$-cocycle with the coefficient in the adjoint representation if and only if $f$ is a derivation on $\g$.
  \end{pro}
  \pf By straightforward computation, for $x_1, \cdots   x_{n}\in\g$, we have
 \begin{eqnarray*}
  &&\delta(f)(x_1\wedge\cdots \wedge x_{n-1},x_n)\\
  &=&-\alpha f([\alpha^{-2}(x_1),\cdots , \alpha^{-2}(x_{n-1}),\alpha^{-2}(x_n)]_\g)+[x_1,\cdots,x_{n-1},f(\alpha^{-1}x_n)]_\g\\
  &&+\sum_{i=1}^{n-1}(-1)^{n-i}[x_1,\cdots,x_{i-1},\hat{x_i},x_{i+1},\cdots,x_{n},f(\alpha^{-1}(x_i))]_\g.
 \end{eqnarray*}
 Thus, $\delta(f)=0$ if and only if
 \begin{eqnarray*}
 && f([\alpha^{-2}(x_1),\cdots , \alpha^{-2}(x_{n-1}),\alpha^{-2}(x_n)]_\g)\\
 &=&
 \sum_{i=1}^{n}[\alpha^{-1} (x_1),\cdots,\alpha^{-1} (x_{i-1}),\alpha^{-1} (f(\alpha^{-1}(x_i))),\alpha^{-1} (x_{i+1}),\cdots,\alpha^{-1} (x_{n})]_\g.
 \end{eqnarray*}
Thus, $f$ is a derivation. The proof is finished. \qed

  \begin{rmk}
    Recall that a derivation of a Lie algebra is exactly a $1$-cocycle with the coefficient in the adjoint representation. Thus, the above proposition justifies our definition of a derivation of an $n$-Hom-Lie algebra.
  \end{rmk}

For all $Y\in\wedge ^{n-1}\g$, $\ad_{Y}$ is a derivation of the $n$-Hom-Lie algebra $(\g,[\cdot,\cdots,\cdot]_\g,\alpha)$, which we call an {\bf inner derivation}. This follows from
\begin{eqnarray*}
\ad_{Y}[x_1,\cdots,x_n]_\g&=&[y_1,\cdots,y_{n-1},[x_1,\cdots,x_n]_\g]_\g\\
&=&[\alpha(\alpha^{-1}y_1),\cdots,\alpha(\alpha^{-1}y_{n-1}),[x_1,\cdots,x_n]_\g]_\g\\
&=&\sum_{i=1}^n[\alpha(x_1),\cdots,\alpha(x_{i-1}),
     [\alpha^{-1}y_1,\cdots,\alpha^{-1}y_{n-1},x_i]_\g,\alpha(x_{i+1}),\cdots,\alpha(x_n)]_\g\\
&=&\sum_{i=1}^n[\alpha(x_1),\cdots,\alpha(x_{i-1}),
     (\alpha^{-1}\circ\ad_{Y}\circ\alpha)(x_i),\alpha(x_{i+1}),\cdots,\alpha(x_n)]_\g.
\end{eqnarray*}
 Denote by $\Inn(\g)$ the set of inner derivations of the $n$-Hom-Lie algebra $(\g,[\cdot,\cdots,\cdot]_\g,\alpha)$, i.e.
\begin{eqnarray}
\Inn(\mathfrak{g})=\{\ad_{Y}\mid Y\in\wedge ^{n-1}\g\}.
\end{eqnarray}

\begin{lem}\label{lem:ideal}
Let $(\mathfrak{g},[\cdot,\cdots,\cdot]_{\mathfrak{g}},\alpha)$ be an $n$-Hom-Lie algebra.  For all $X=x_1\wedge\cdots\wedge x_{n-1}\in\wedge^{n-1}\g$ and $D\in\Der(\g)$, we have
  $$\Ad_{\alpha}\ad_X=\ad_{\widetilde{\alpha}(X)},\quad [D,\ad_X]_{\alpha}=\ad_{\sum_{i=1}^{n-1}\alpha(x_1)\wedge\cdots\wedge D(x_{i})\wedge\cdots\alpha(x_{n-1})}.$$
  Therefore, $\Inn(\g)$ is an ideal of the Hom-Lie algebra $(\Der(\g),[\cdot,\cdot]_{\alpha},\Ad_{\alpha})$.
\end{lem}
\pf For all $y\in \g$, we have 
\begin{eqnarray*}
(\Ad_{\alpha}\ad_X)(y)&=&(\alpha \circ \ad_X \circ \alpha^{-1})(y)
                       =\alpha[x_1,\cdots,x_{n-1},\alpha^{-1}(y)]_\g
                       =[\alpha(x_1),\cdots,\alpha(x_{n-1}),y]_\g\\
                      &=&\ad_{\widetilde{\alpha}(X)}(y).
\end{eqnarray*}
By \eqref{defi:derivation}, we have
\begin{eqnarray*}
&&[D,\ad_{X}]_{\alpha}(y)\\
                &=&(\alpha\circ D\circ\alpha^{-1}\circ\ad_{X}\circ\alpha^{-1})(y)-
                (\alpha\circ\ad_{X}\circ\alpha^{-1}\circ D\circ\alpha^{-1})(y)\\
                &=&\alpha(D[\alpha^{-1}(x_1),\cdots,\alpha^{-1}(x_{n-1}),\alpha^{-2}(y)]_{\mathfrak{g}})-
                \alpha[x_1,\cdots,x_{n-1},\alpha^{-1}D(\alpha^{-1}(y))]_{\mathfrak{g}}\\
                &=&\alpha\Big([x_1,\cdots,x_{n-1},(\Ad_{\alpha^{-1}}D)(\alpha^{-2}(y))]_{\mathfrak{g}}
                +\sum_{i=1}^{n-1}[x_1,\cdots,(\Ad_{\alpha^{-1}}D)(\alpha^{-1}(x_i)),\cdots,x_{n-1},\alpha^{-1}(y)]_{\mathfrak{g}}\Big)
                \\
                &&-[\alpha(x_1),\cdots,\alpha(x_{n-1}),D(\alpha^{-1}(y))]_{\mathfrak{g}}\\
                                &=&\ad_{\sum_{i=1}^{n-1}\alpha(x_1)\wedge\cdots\wedge D(x_{i})\wedge\cdots\alpha(x_{n-1})}y.
\end{eqnarray*}
 The proof is finished.\qed

 \subsection{Dual representations of $n$-Hom-Lie algebras}

 In this subsection, we study dual representation of an $n$-Hom-Lie algebra. See \cite{caisheng} for more details about dual representations of Hom-Lie algebras.

 Let $(V,\rho,\beta)$  be a representation of an $n$-Hom-Lie algebra $(\g,[\cdot,\cdots,\cdot]_\g,\alpha)$, where $\beta$ is invertible, i.e. $\beta\in GL(V)$.
Define $\rho^*:\wedge^{n-1}\g\rightarrow\mathfrak{gl}(V^*)$ as usual by
\begin{equation}\label{eq:defirhodual}
\langle\rho^*(X)(\xi),u\rangle=-\langle\xi,\rho(X)(u)\rangle,\quad\forall X\in \wedge^{n-1}\g,~u\in V,~\xi\in V^*.
\end{equation}
Then we define $\rho^\star:\wedge^{n-1}\g\rightarrow\gl(V^*)$ by
\begin{equation}\label{eq:defirhodualnew}
\rho^\star(X)(\xi)=\rho^*(\widetilde{\alpha}(X))((\beta^{-2})^*(\xi)).
\end{equation}

\begin{thm}
Let $(V,\rho,\beta)$  be a representation of an $n$-Hom-Lie algebra $(\g,[\cdot,\cdots,\cdot]_\g,\alpha)$, where $\beta$ is invertible. Then $(V^*,\rho^\star,(\beta^{-1})^*)$ is also a representation of $(\g,[\cdot,\cdots,\cdot]_\g,\alpha)$, which we call the {\bf dual representation} of the representation $(V,\rho,\beta)$.
\end{thm}
\pf By $\rho(\widetilde{\alpha}(X) )\circ\beta=\beta\circ \rho(X)$, we have
\begin{eqnarray*}
\langle\rho^\star(\widetilde{\alpha}(X))((\beta^{-1})^*(\xi)),u\rangle
&=&\langle\rho^*(\widetilde{\alpha}^2(X))((\beta^{-3})^*(\xi)),u\rangle=-\langle(\beta^{-3})^*(\xi),\rho(\widetilde{\alpha}^2(X))(u)\rangle\\
&=&-\langle(\beta^{-2})^*(\xi),\beta^{-1}\rho(\widetilde{\alpha}^2(X))(u)\rangle=-\langle(\beta^{-2})^*(\xi),\rho(\widetilde{\alpha}(X))(\beta^{-1}(u))\rangle\\
&=&\langle\rho^*(\widetilde{\alpha}(X))((\beta^{-2})^*(\xi)),\beta^{-1}(u)\rangle=\langle\rho^\star(X)(\xi),\beta^{-1}(u)\rangle\\
&=&\langle(\beta^{-1})^*\rho^\star(X)(\xi),u\rangle,
\end{eqnarray*}
which implies that
 $$\rho^\star(\widetilde{\alpha}(X))\circ(\beta^{-1})^*=(\beta^{-1})^*\circ\rho^\star(X).$$

By straightforward computation, we have
\begin{eqnarray*}
\langle\rho^\star(\widetilde{\alpha}(X))\big{(}\rho^\star(Y)(\xi)\big{)},u\rangle
&=&\langle\rho^*(\widetilde{\alpha}^2(X))\big{(}(\beta^{-2})^*
\rho^*(\widetilde{\alpha}(Y))((\beta^{-2})^*(\xi))\big{)},u\rangle\\
&=&-\langle(\beta^{-2})^*\rho^*(\widetilde{\alpha}(Y))((\beta^{-2})^*(\xi)),\rho(\widetilde{\alpha}^2(X))(u)\rangle\\
&=&-\langle\rho^*(\widetilde{\alpha}(Y))((\beta^{-2})^*(\xi)),\beta^{-2}\rho(\widetilde{\alpha}^2(X))(u)\rangle\\
&=&\langle(\beta^{-2})^*(\xi),\rho(\widetilde{\alpha}(Y))\big{(}\beta^{-2}\rho(\widetilde{\alpha}^2(X))(u)\big{)}\rangle\\
&=&\langle(\beta^{-2})^*(\xi),\rho(\widetilde{\alpha}(Y))\big{(}\rho(X)(\beta^{-2}(u))\big{)}\rangle.
\end{eqnarray*}
Therefore, we have
\begin{eqnarray*}
&&\langle\big{(}\rho^\star(\widetilde{\alpha}(X))\circ\rho^\star(Y)
-\rho^\star(\widetilde{\alpha}(Y))\circ\rho^\star(X)\big{)}(\xi),u\rangle\\
&=&\langle(\beta^{-2})^*(\xi),\rho(\widetilde{\alpha}(Y))\big{(}\rho(X)(\beta^{-2}(u))\big{)}
-\rho(\widetilde{\alpha}(X))\big{(}\rho(Y)(\beta^{-2}(u))\big{)}\rangle\\
&=&\langle(\beta^{-2})^*(\xi),-\big{(}\rho(\widetilde{\alpha}(X))\circ\rho(Y)
-\rho(\widetilde{\alpha}(Y))\circ\rho(X)\big{)}(\beta^{-2}(u))\rangle\\
&=&\langle(\beta^{-2})^*(\xi),
-\sum_{i=1}^{n-1}\rho(\alpha(y_1),\cdots,\alpha(y_{i-1}),\ad_X(y_i),\alpha(y_{i+1}),
\cdots,\alpha(y_{n-1}))(\beta^{-1}(u))\rangle\\
&=&-\sum_{i=1}^{n-1}\langle(\beta^{-3})^*(\xi),
\beta\rho(\alpha(y_1),\cdots,\alpha(y_{i-1}),\ad_X(y_i),\alpha(y_{i+1}),
\cdots,\alpha(y_{n-1}))(\beta^{-1}(u))\rangle\\
&=&-\sum_{i=1}^{n-1}\langle(\beta^{-3})^*(\xi),
\rho(\alpha^2(y_1),\cdots,\alpha^2(y_{i-1}),\alpha(\ad_X(y_i)),\alpha^2(y_{i+1}),
\cdots,\alpha^2(y_{n-1}))(u)\rangle\\
&=&\sum_{i=1}^{n-1}\langle\rho^*(\alpha^2(y_1),\cdots,\alpha^2(y_{i-1}),\alpha(\ad_X(y_i)),\alpha^2(y_{i+1}),
\cdots,\alpha^2(y_{n-1}))((\beta^{-3})^*(\xi)),u\rangle\\
&=&\langle\sum_{i=1}^{n-1}\rho^\star(\alpha(y_1),\cdots,\alpha(y_{i-1}),
\ad_X(y_i),\alpha(y_{i+1}),\cdots,\alpha(y_{n-1}))((\beta^{-1})^*(\xi)),u\rangle,
\end{eqnarray*}
which implies that
\begin{eqnarray*}
\rho^\star(\widetilde{\alpha}(X))\circ\rho^\star(Y)-\rho^\star(\widetilde{\alpha}(Y))\circ\rho^\star(X)=\sum_{i=1}^{n-1}\rho^\star([X,Y]_F)\circ(\beta^{-1})^*.
\end{eqnarray*}

Finally, by straightforward computations, we can obtain
\begin{eqnarray*}
        &&\rho^\star(\alpha(x_1),\cdots,\alpha(x_{n-2}),[y_1,\cdots,y_n]_\g)\circ(\beta^{-1})^*\\
        &=&\sum_{i=1}^n(-1)^{n-i}\rho^\star(\alpha(y_1),\cdots,\hat{y_i},\cdots,\alpha(y_n))\circ\rho^\star(x_1,\cdots,x_{n-2},y_i).
        \end{eqnarray*}

Therefore, $(V^*,\rho^\star,(\beta^{-1})^*)$ is a representation of $(\g,[\cdot,\cdots,\cdot]_\g,\alpha)$. \qed

\begin{rmk}
  By straightforward computation, one can deduce  that the usual $\rho^*$ defined by \eqref{eq:defirhodual} is not a representation anymore. In particular, for the adjoint representation $\ad$, the dual map $\ad^*$ is not a representation in general.  On the other hand, to solve this problem, the authors in \cite{bai} add some strong conditions to make $\rho^*$ still being a representation. Here, we see that there is a natural definition of $\rho^\star$ such that it is a representation. This makes our definition of $\rho^\star$ nontrivial. Furthermore, the definition of $\rho^\star$ is a generalization of the one given in \cite{caisheng} for Hom-Lie algebras.
\end{rmk}

   \section{$(n-1)$-order deformations, Hom-Nijenhuis operators and Hom-$\huaO$-operators of an $n$-Hom-Lie algebra}

\subsection{$(n-1)$-order deformations of an $n$-Hom-Lie algebra}

Let $(\g,[\cdot,\cdots,\cdot]_\g,\alpha)$ be an $n$-Hom-Lie algebra. For convenience, we denote by $\omega_0(x_1,\cdots,x_n)=[x_1,\cdots,x_n]_\g$. Let $\omega_i:\wedge^{n}\g\longrightarrow \g,~~ 1\leq i\leq n-1$ be skew-symmetric   multilinear   maps. Consider a $\lambda$-parametrized family of linear operations:

\begin{equation}\label{eq:deformation}
\omega_\lambda(x_1,\cdots,x_n)=\sum_{i=0}^{n-1}\lambda^i\omega_i(x_1,\cdots,x_n).
\end{equation}
Here $\lambda\in\mathbb K,$ where $\mathbb K$ is the base field. If all $(\g,\omega_\lambda,\alpha)$ are  $n$-Hom-Lie algebras, we say that
$\omega_1,\cdots,\omega_{n-1}$ generate   {\bf an $(n-1)$-order $1$-parameter  deformation} of the $n$-Hom-Lie algebra $(\g,\omega_0,\alpha)$.
We also denote by $[x_1,\cdots,x_n]_\lambda=\omega_\lambda(x_1,\cdots,x_n)$.

\begin{pro}\label{conds}
With the above notations, $\omega_1,\cdots,\omega_{n-1}$ generate   an $(n-1)$-order $1$-parameter   deformation of the $n$-Hom-Lie algebra $(\g,\omega_0,\alpha)$ if and only if for all $i,j=1,2,\cdots,n-1$ and $k=1,2,\cdots,2n-2$ the following conditions are satisfied:
  \begin{eqnarray}
\label{hom-n-lie-1}\omega_i\circ \alpha^{\otimes n}&=&\alpha\circ\omega_i,\\
\label{hom-n-lie-2}\sum\limits_{i+j=k}\omega_i\circ\omega_j&=&0.
\end{eqnarray}
Here $\omega_i\circ\omega_j:\wedge^{n-1}\g\otimes \wedge^{n-1}\g\wedge\g\longrightarrow\g$ is defined by
\begin{eqnarray}\label{N-R bracket}
\nonumber\omega_i\circ\omega_j(X,Y,z)=\omega_i(\omega_j(X,\cdot)*Y,\alpha(z))
-\omega_i(\widetilde{\alpha}(X),\omega_j(Y,z))+\omega_i(\widetilde{\alpha}(Y),\omega_j(X,z)),
\end{eqnarray}
where $\omega_j(X,\cdot)*Y\in\wedge^{n-1}\g$ is given by
\begin{eqnarray}
\omega_j(X,\cdot)*Y=\sum_{l=1}^{n-1}\alpha(y_1)\wedge\cdots\wedge\alpha(y_{l-1})\wedge\omega_j(X,y_l)\wedge\alpha(y_{l+1})\wedge\cdots\wedge\alpha(y_{n-1}).
\end{eqnarray}
\end{pro}
\pf
$(\g,\omega_\lambda,\alpha)$ are $n$-Hom-Lie algebra structures if and only if
\begin{eqnarray}
\label{deformation-1}\omega_\lambda\circ \alpha^{\otimes n}&=&\alpha\circ\omega_\lambda,\\
\label{deformation-2}\omega_\lambda(\widetilde{\alpha}(X),\omega_\lambda(Y,z))&=&\omega_\lambda(\omega_\lambda(X,\cdot)*Y,\alpha(z))+\omega_\lambda(\widetilde{\alpha}(Y),\omega_\lambda(X,z)).
\end{eqnarray}
By \eqref{deformation-1}, we have
$$\omega_i\circ \alpha^{\otimes n}=\alpha\circ\omega_i.$$
Expanding the equations in \eqref{deformation-2} and collecting coefficients of $\lambda^k$, we see that  \eqref{deformation-2} is equivalent to the system of equations
\begin{eqnarray*}
\sum_{i+j=k}\omega_i(\widetilde{\alpha}(X),\omega_j(Y,z))=\sum_{i+j=k}\omega_i(\omega_j(X,\cdot)*Y,\alpha(z))+\sum_{i+j=k}\omega_i(\widetilde{\alpha}(Y),\omega_j(X,z)).
\end{eqnarray*}
Thus, we have
\begin{eqnarray*}
\sum_{i+j=k}\omega_i(\omega_j(X,\cdot)*Y,\alpha(z))-\omega_i(\widetilde{\alpha}(X),\omega_j(Y,z))+\omega_i(\widetilde{\alpha}(Y),\omega_j(X,z)),
\end{eqnarray*}
which finishes the proof. \qed

\begin{cor}
  If  $\omega_1,\cdots,\omega_{n-1}$ generate   an $(n-1)$-order $1$-parameter   deformation of the $n$-Hom-Lie algebra $(\g,\omega_0,\alpha)$, then $\omega_1$ is a $2$-cocycle of the $n$-Hom-Lie algebra $(\g,\omega_0,\alpha)$ with the coefficients in the adjoint representation.
\end{cor}
\pf By \eqref{hom-n-lie-2}, let $k$=1, we deduce that
$$
\omega_0\circ\omega_1+\omega_1\circ\omega_0=0,
$$
which is equivalent to that $\omega_1$ is a 2-cocycle. We omit details. \qed

\begin{cor}
  If  $\omega_1,\cdots,\omega_{n-1}$ generate   an $(n-1)$-order $1$-parameter   deformation of the $n$-Hom-Lie algebra $(\g,\omega_0,\alpha)$, then  $(\g,\omega_{n-1},\alpha)$ is an $n$-Hom-Lie algebra.
\end{cor}
\pf By \eqref{hom-n-lie-1}, let $i=n-1$, we deduce that
$$\omega_{n-1}\circ \alpha^{\otimes n}=\alpha\circ\omega_{n-1},$$
and by \eqref{hom-n-lie-2}, let $k=2n-2$, we deduce that
\begin{eqnarray*}
\omega_{n-1}\circ\omega_{n-1}=0,
\end{eqnarray*}
which is equivalent to that $(\g,\omega_{n-1},\alpha)$ is an $n$-Hom-Lie algebra.  \qed

\subsection{Hom-Nijenhuis operators and Hom-$\huaO$-operators of an $n$-Hom-Lie algebra}

 In this subsection,  we study trivial $(n-1)$-order deformations of an $n$-Hom-Lie algebra and  introduce the notion of a Hom-Nijenhuis operator of an $n$-Hom-Lie algebra, which could  generate a trivial deformation.
 Then we give the relation between Hom-$\mathcal O$-operators and Hom-Nijenhuis operators.

 \begin{defi}\label{defi:trivial}
A deformation is said to be {\bf trivial} if there exists a linear map $N:\g\longrightarrow \g$ such that for all $\lambda$,  $T_\lambda=\alpha+\lambda N$ satisfies
\begin{eqnarray}
\label{trivial-1}T_\lambda\circ\alpha&=&\alpha\circ T_\lambda,\\
\label{trivial-2}T_\lambda[x_1, \cdots,x_n]_\lambda&=&[ T_\lambda x_1,  \cdots, T_\lambda x_n]_\g,\quad \forall x_1,\cdots, x_n\in \g.
\end{eqnarray}
\end{defi}

Eq.~\eqref{trivial-1} equals to
\begin{eqnarray}
N\circ\alpha&=&\alpha\circ N.
\end{eqnarray}
The left hand side of Eq.~\eqref{trivial-2} equals to
\begin{eqnarray*}
&&\alpha[x_1, \cdots,x_n]_\g+\lambda\big(\alpha\omega_1(x_1, \cdots,x_n)+N[x_1, \cdots,x_n]_\g\big)\\
&&+\sum_{i=2}^{n-1}\lambda^{i}\big(\alpha\omega_{i}(x_1,\cdots,x_n)+N\omega_{i-1}(x_1,\cdots,x_n)\big)+\lambda^nN\omega_{n-1}(x_1,\cdots,x_n).
\end{eqnarray*}
The right hand side of Eq.~\eqref{trivial-2} equals to
\begin{eqnarray*}
&&[\alpha(x_1),\cdots,\alpha(x_n)]_\g+\lambda\sum_{i=1}^n[\alpha(x_1),\cdots,Nx_i,\cdots,\alpha(x_n)]_\g\\
&&
+\lambda^2\sum_{l_1<l_2}[\alpha(x_1),\cdots,Nx_{l_1},\cdots,Nx_{l_2},\cdots,\alpha(x_n)]_\g+\cdots+\lambda^n[Nx_1,\cdots,Nx_n]_\g.
\end{eqnarray*}
Therefore, by Eq.~\eqref{trivial-2}, we have
\begin{eqnarray}\label{eq:trivialdef con1}
 \alpha\omega_1(x_1,\cdots,x_n)+N[x_1,\cdots,x_n]_\g&=&\sum_{i=1}^n[\alpha(x_1),\cdots,Nx_i,\cdots,\alpha(x_n)]_\g,\\
 \label{eq:trivialdef con2}N\omega_{n-1}(x_1,\cdots,x_n)&=&[Nx_1,\cdots,Nx_n]_\g,
 \end{eqnarray}
and
\begin{eqnarray}\label{eq:trivialdef con3}
\nonumber&&\alpha \omega_i(x_1,\cdots,x_n)+N\omega_{i-1}(x_1,\cdots,x_n)\\
&=&\sum_{l_1<l_2\cdots<l_i}[\alpha(x_1),\cdots,Nx_{l_1},\cdots,Nx_{l_2},\cdots,Nx_{l_i},\cdots,\alpha(x_n)]_\g,
\end{eqnarray}
for all $2\leq i\leq n-1$.

\emptycomment{
\begin{eqnarray*}\label{eq:trivialdef con1}
 \omega_1(x_1,x_2,\cdots,x_n)+N[x_1,x_2,\cdots,x_n]&=&\sum_{i=1}^n[x_1,\cdots,Nx_i,\cdots,x_n],\\
 \omega_2(x_1,x_2,\cdots,x_n)+N\omega_1(x_1,x_2,\cdots,x_n)&=&\sum_{i<j}[x_1,\cdots,Nx_i,\cdots,Nx_j,\cdots,x_n],\\
 \omega_3(x_1,x_2,\cdots,x_n)+N\omega_2(x_1,x_2,\cdots,x_n)&=&\sum_{i<j<k}[x_1,\cdots,Nx_i,\cdots,Nx_j,\cdots,Nx_k,\cdots,x_n],\\
  \quad\quad\quad\quad\quad\quad\quad\quad\quad\quad\quad\quad\vdots  \\
 \omega_l(x_1,x_2,\cdots,x_n)+N\omega_{l-1}(x_1,x_2,\cdots,x_n)
&=&\sum_{i_1<i_2\cdots<i_l}[x_1,\cdots,Nx_{i_1},\cdots,Nx_{i_k},\cdots,Nx_{i_l},\cdots,x_n],\\
  \quad\quad\quad\quad\quad\quad\quad\quad\quad\quad\quad\quad\vdots  \\
   N\omega_{n-1}(x_1,x_2,\cdots,x_n)&=&[Nx_1,Nx_2,\cdots,Nx_n].
\end{eqnarray*}
}

Let $(\g,[\cdot,\cdots,\cdot]_\g,\alpha)$ be an $n$-Hom-Lie algebra, and $N:\g\longrightarrow\g$ a linear map. Define an $n$-ary bracket $[\cdot,\cdots,\cdot]_N^1:\wedge^n\g\longrightarrow\g$  by
   \begin{equation}\label{eq:bracket(1)}
    [x_1,\cdots,x_n]_N^{1}=\sum_{i=1}^n[x_1,\cdots,N\alpha^{-1}(x_i),\cdots,x_n]_\g-N[\alpha^{-1}(x_1),\cdots,\alpha^{-1}(x_n)]_\g.
  \end{equation}
  Then we define $n$-ary brackets $[\cdot,\cdots,\cdot]_N^i:\wedge^n\g\longrightarrow\g, (2\leq i\leq n-1)$ via induction by
    \begin{eqnarray}\label{eq:bracket (j)}
   \nonumber [x_1,\cdots,x_n]_N^{i}&=&
    \sum_{l_1<l_2\cdots<l_{i}}[x_1,\cdots,N\alpha^{-1}(x_{l_1}),\cdots,N\alpha^{-1}(x_{l_2}),\cdots,N\alpha^{-1}(x_{l_i}),\cdots,x_n]_\g\\
    &&-N[\alpha^{-1}(x_1), \cdots,\alpha^{-1}(x_n)]_N^{i-1}.
  \end{eqnarray}

\begin{defi}\label{defi:Nijenhuis}
Let $(\g,[\cdot,\cdots,\cdot]_\g,\alpha)$ be an $n$-Hom-Lie algebra. A linear map    $N:\g\longrightarrow\g$   is called a {\bf Hom-Nijenhuis operator} if
\begin{eqnarray}
\label{eq:Nijenhuis(n)-1}N\circ\alpha              &=&\alpha\circ N,\\
\label{eq:Nijenhuis(n)-2}[Nx_1,\cdots,Nx_n]_\g&=&N( [x_1,\cdots,x_n]_N^{n-1}),\quad\forall x_1,\cdots, x_n\in \g.
\end{eqnarray}
\end{defi}

\begin{rmk}
  When $n=2$,  \eqref{eq:Nijenhuis(n)-2} reduces to
  $$
  [Nx,Ny]_\g=N[x,y]_N^1=N([x,N\alpha^{-1}(y)]_\g+[N\alpha^{-1}(x),y]_\g-N[\alpha^{-1}(x),\alpha^{-1}(y)]_\g),
  $$
  which is the same as the condition given in \cite[Proposition 6.2]{CaiSheng}. Thus, the notion of a Hom-Nijenhuis operator of an $n$-Hom-Lie algebra is a natural generalization of the Hom-Nijenhuis operator of a Hom-Lie algebra given in \cite{CaiSheng}.
\end{rmk}

\begin{thm}\label{thm:trivial deformation}
  Let $N$ be a Hom-Nijenhuis operator of an $n$-Hom-Lie algebra
$(\g,[\cdot,\cdots,\cdot]_\g,\alpha)$. Then
   a deformation can be obtained by putting
  \begin{equation}\label{eq:formulaOT}
    \omega_i(x_1,x_2,\cdots,x_n)=[x_1,x_2,\cdots,x_n]_N^{i},\quad 1\leq i\leq n-1.
  \end{equation}
  Moreover, this deformation  is trivial.
\end{thm}

\begin{lem}\label{lem:isomorphism}
  Let $(\g,[\cdot,\cdots,\cdot]_\g,\alpha)$ be an $n$-Hom-Lie algebra and $\h$ a vector space with a linear automorphism $\gamma\in GL(\h)$. If there exists an linear isomorphism $f:\h\longrightarrow \g$, such that $\alpha\circ f=f\circ\gamma $. We define an $n$-ary bracket $[\cdot,\cdots,\cdot]'$ on $\h$ by
  $$
 [u_1,u_2,\cdots,u_n]'= f^{-1}[f(u_1),f(u_2),\cdots,f(u_n)]_\g,\quad\forall ~u_i\in\h,
$$
then $(\h,[\cdot,\cdots,\cdot]',\gamma)$ is an $n$-Hom-Lie algebra.
\end{lem}
\pf It follows from
straightforward computations.   \qed\vspace{3mm}

{\bf The proof of Theorem \ref{thm:trivial deformation}:}
It is obvious that for a Hom-Nijenhuis   operator $N$,
the maps $\omega_1,\cdots, \omega_{n-1}$ given by Eq.~\eqref{eq:formulaOT} satisfy Eq.~\eqref{eq:trivialdef con1}-\eqref{eq:trivialdef con3}. Therefore,
for any $\lambda$, $T_\lambda$ satisfies
\begin{eqnarray*}
\label{eq:trivial-1}T_\lambda\circ\alpha&=&\alpha\circ T_\lambda,\\
\label{eq:trivial-2}T_\lambda[x_1,x_2,\cdots,x_n]_\lambda&=&[ T_\lambda x_1, T_\lambda x_2, \cdots, T_\lambda x_n]_\g,\quad \forall x_1,\cdots, x_n\in \g.
\end{eqnarray*}
 For $\lambda$ sufficiently small, we see that $T_\lambda$ is an
isomorphism between vector spaces. Thus, we have
$$[x_1,x_2,\cdots,x_n]_\lambda=T_\lambda^{-1}[ T_\lambda x_1, T_\lambda x_2, \cdots, T_\lambda x_n]_\g.$$
By Lemma \ref{lem:isomorphism},
we deduce that $(\g,[\cdot,\cdots,\cdot]_\lambda,\alpha)$ is an   $n$-Hom-Lie algebra, for $\lambda$ sufficiently small.
  Thus, $\omega_1,\cdots, \omega_{n-1}$
  given by
Eq.~\eqref{eq:formulaOT} satisfy the conditions
(\ref{hom-n-lie-1})-(\ref{hom-n-lie-2}) in Proposition \ref{conds}. Therefore,
$(\g,[\cdot,\cdots,\cdot]_\lambda,\alpha)$ is an $n$-Hom-Lie algebra for all
$\lambda$, which means that $\omega_1,\cdots, \omega_{n-1}$
given by
Eq.~\eqref{eq:formulaOT} generate a deformation. It is obvious  that
this deformation is trivial.\qed

\begin{cor}\label{cor:n-LA}
  Let $N$ be a Hom-Nijenhuis operator of an $n$-Hom-Lie algebra
$(\g,[\cdot,\cdots,\cdot]_\g, \alpha)$.
  Then $(\g,[\cdot,\cdots,\cdot]_N^{n-1},\alpha)$ is an $n$-Hom-Lie algebra, and $N$ is a   homomorphism from $(\g,[\cdot,\cdots,\cdot]_N^{n-1},\alpha)$ to $(\g,[\cdot,\cdots,\cdot]_\g,\alpha)$.
\end{cor}

In the sequel, we introduce the notion of a Hom-$\huaO$-operator associated to a representation of an $n$-Hom-Lie algebra and give the relation between Hom-$\huaO$-operators and Hom-Nijenhuis operators.

\begin{defi}
Let $(\g,[\cdot,\cdots,\cdot]_\g,\alpha)$ be an $n$-Hom-Lie algebra and $(V,\rho,\beta)$ a representation. A linear map $T:V\lon\g$ is called a {\bf Hom-$\huaO$-operator} if for all $v_1,\cdots,v_n\in V$,
\begin{eqnarray}
\label{o1}T\circ\beta&=&\alpha\circ T,\\
\label{o2}[Tv_1,\cdots,Tv_n]_\g&=&T\sum_{i=1}^{n}(-1)^{n-i}\rho\big(T\beta^{-1}(v_1),\cdots, \widehat{T\beta^{-1}(v_{i})}, \cdots,T\beta^{-1}(v_{n})\big)(v_i).
\end{eqnarray}
\end{defi}

\begin{pro}
Let $(\g,[\cdot,\cdots,\cdot]_\g,\alpha)$ be an $n$-Hom-Lie algebra and $(V,\rho,\beta)$ a representation. A linear map $T:V\lon\g$ is   a  Hom-$\huaO$-operator if and only if
\begin{eqnarray*}
\bar{T}=\left(\begin{array}{cc}0&T\\0&0\end{array}\right):\g\oplus V\lon\g\oplus V
\end{eqnarray*}
is a Hom-Nijenhuis operator acting on the semidirect product $n$-Hom-Lie algebra $\g\ltimes V$.
\end{pro}
\pf First it is obvious that $\bar{T}\circ(\alpha+\beta)=(\alpha+\beta)\circ\bar{T} $ if and only if
$
T\circ\beta=\alpha\circ T.
$
Then for all $x_1,\cdots,x_n\in\g,v_1,\cdots,v_n\in V$, we have
\begin{eqnarray*}
[\bar{T}(x_1+v_1),\bar{T}(x_2+v_2),\cdots,\bar{T}(x_n+v_n)]_{\rho}=[Tv_1,Tv_2,\cdots,Tv_n]_\g.
\end{eqnarray*}
On the other hand, since $\bar{T}^2=0$, we have
\begin{eqnarray*}
&&\bar{T}[x_1+v_1,x_2+v_2,\cdots,x_n+v_n]_{\bar{T}}^{n-1}\\
&=&\bar{T}\big(\sum_{l_1<l_2\cdots<l_{n-1}}[\cdots,\bar{T}(\alpha^{-1}(x_{l_1})+\beta^{-1}(v_{l_1})),\cdots,\bar{T}(\alpha^{-1}(x_{l_{n-1}})
+\beta^{-1}(v_{l_{n-1}})),\cdots]_{\rho}\big)\\
&=&\bar{T}\sum_{i=1}^{n}[T\beta^{-1}(v_1),\cdots,T\beta^{-1}(v_{i-1}),x_i+v_i,T\beta^{-1}(v_{i+1}),\cdots,T\beta^{-1}(v_{n})]_{\rho}\\
&=&T\sum_{i=1}^{n}(-1)^{n-i}\rho(T\beta^{-1}(v_1),\cdots,T\beta^{-1}(v_{i-1}),\widehat{T\beta^{-1}(v_{i})},T\beta^{-1}(v_{i+1}),\cdots,T\beta^{-1}(v_{n}))(v_i),
\end{eqnarray*}
which implies that $\bar{T}$ is a Hom-Nijenhuis operator if and only if  $T$ is a Hom-$\huaO$-operator. \qed

 \section{Generalized derivation extensions of $n$-Hom-Lie algebras}

In this section, we give a new approach to construct $n$-Hom-Lie algebras, which is called  generalized derivation extensions of $n$-Hom-Lie algebras. First we give the notion of a generalized derivation of an $n$-Hom-Lie algebra.

\begin{defi}\label{derivation}
Let $(\g,[\cdot,\cdots,\cdot]_\g,\alpha)$ be an $n$-Hom-Lie algebra. A linear map $D:\wedge^{n-1}\g\lon\g$ is called a {\bf generalized derivation}, if for all $x_1,\cdots,x_{n-1}, y_1,\cdots,y_{n-1}\in\g$, the following conditions are satisfied:
\begin{itemize}
  \item[{\rm(i)}] $\alpha\circ D=D\circ\widetilde{\alpha};$

  \item[{\rm(ii)}]
  \begin{eqnarray*}
  &&[\alpha(x_1),\cdots,\alpha(x_{n-1}),D(y_1,\cdots,y_{n-1})]_\g-[\alpha(y_{1}),\cdots,\alpha(y_{n-1}),D(x_1,\cdots,x_{n-1})]_\g\\
  &=&\sum_{i=1}^{n-1}D(\alpha(y_1),\cdots,\alpha(y_{i-1}),
     [x_1,\cdots,x_{n-1},y_i]_\g,\alpha(y_{i+1}),\cdots,\alpha(y_{n-1}));
  \end{eqnarray*}
    \item[{\rm(iii)}] \begin{eqnarray*}
        &&D(\alpha(x_1),\cdots,\alpha(x_{n-2}),[y_1,\cdots,y_{n}]_\g)\\
        &=&\sum_{i=1}^{n}[\alpha(y_1),\cdots,\alpha(y_{i-1}),
     D(x_1,\cdots,x_{n-2},y_i),\alpha(y_{i+1}),\cdots,\alpha(y_{n})]_\g;
        \end{eqnarray*}
  \item[{\rm(iv)}] \begin{eqnarray*}
        &&D(\alpha(x_1),\cdots,\alpha(x_{n-2}),D(y_1,\cdots,y_{n-1}))\\
        &=&\sum_{i=1}^{n-1}D(\alpha(y_1),\cdots,\alpha(y_{i-1}),
     D(x_1,\cdots,x_{n-2},y_i),\alpha(y_{i+1}),\cdots,\alpha(y_{n-1})).
        \end{eqnarray*}
\end{itemize}
\end{defi}

\begin{rmk}
Let $D:\wedge^{n-1}\g\lon\g$ be a generalized derivation on  $(\g,[\cdot,\cdots,\cdot]_\g,\alpha)$. By Conditions (i) and (iv) in Definition \ref{derivation},   $D$ defines an $(n-1)$-Hom-Lie algebra structure on the vector space $\g$.
\end{rmk}

For all $x\in\g$, satisfying $\alpha(x)=x$, define $\add_x:\wedge^{n-1}\g\lon\g$ by
$$\add_x(y_1,\cdots,y_{n-1})=[x,y_1,\cdots,y_{n-1}]_\g.$$
Then we have

\begin{lem}
For all $x\in\g$, satisfying $\alpha(x)=x$, $\add_x$ is a generalized derivation on the $n$-Hom-Lie algebra $(\g,[\cdot,\cdots,\cdot]_\g,\alpha)$, which is called an inner generalized derivation.
\end{lem}

\pf   By $\alpha(x)=x$ we have
\begin{eqnarray*}
(\alpha\circ\add_{x})(y_1,\cdots,y_{n-1})&=&[\alpha(x),\alpha(y_1),\cdots,\alpha(y_{n-1})]_\g\\
                                   &=&[x,\alpha(y_1),\cdots,\alpha(y_{n-1})]_\g\\
                                   &=&(\add_x\circ\widetilde{\alpha})(y_1,\cdots,y_{n-1}).
\end{eqnarray*}
Thus, Condition (i) in Definition \ref{derivation} holds.

By $\alpha(x)=x$ and the Hom Fundamental identity, we have
\begin{eqnarray*}
&&[\alpha(x_1),\cdots,\alpha(x_{n-1}),\add_x(y_1,\cdots,y_{n-1})]_\g-[\alpha(y_1),\cdots,\alpha(y_{n-1}),\add_x(x_1,\cdots,x_{n-1})]_\g\\
&&=[\alpha(x_1),\cdots,\alpha(x_{n-1}),[x,y_1,\cdots,y_{n-1}]_\g]_\g-[\alpha(y_1),\cdots,\alpha(y_{n-1}),[x,x_1,\cdots,x_{n-1}]_\g]_\g\\
&&=\sum_{i=1}^{n-1}[\alpha(x),\alpha(y_1),\cdots,\alpha(y_{i-1}),[x_1,\cdots,x_{n-1},y_i]_\g,\alpha(y_{i+1}),\cdots,\alpha(y_{n-1})]\\
&&=\sum_{i=1}^{n-1}\add_x(\alpha(y_1),\cdots,\alpha(y_{i-1}),[x_1,\cdots,x_{n-1},y_i]_\g,\alpha(y_{i+1}),\cdots,\alpha(y_{n-1})).
\end{eqnarray*}
Thus, Condition (ii) in Definition \ref{derivation} holds.

Similarly,   we have
\begin{eqnarray*}
&&\add_x(\alpha(x_1),\cdots,\alpha(x_{n-2}),[y_1,\cdots,y_{n}]_\g)\\
&=&[\alpha(x),\alpha(x_1),\cdots,\alpha(x_{n-2}),[y_1,\cdots,y_{n}]_\g]_\g\\
&=&\sum_{i=1}^{n}[\alpha(y_1),\cdots,\alpha(y_{i-1}),[x,x_1,\cdots,x_{n-2},y_i]_\g,\alpha(y_{i+1}),\cdots,\alpha(y_{n})]_\g\\
&=&\sum_{i=1}^{n}[\alpha(y_1),\cdots,\alpha(y_{i-1}),\add_x(x_1,\cdots,x_{n-2},y_i),\alpha(y_{i+1}),\cdots,\alpha(y_{n})]_\g,
\end{eqnarray*}
 and
\begin{eqnarray*}
&&\add_x(\alpha(x_1),\cdots,\alpha(x_{n-2}),\add_x(y_1,\cdots,y_{n-1}))\\
&=&[\alpha(x),\alpha(x_1),\cdots,\alpha(x_{n-2}),[x,y_1,\cdots,y_{n-1}]_\g]_\g\\
&=&\sum_{i=1}^{n-1}[x,\alpha(y_1),\cdots,\alpha(y_{i-1}),[x,x_1,\cdots,x_{n-2},y_i]_\g,\alpha(y_{i+1}),\cdots,\alpha(y_{n})]_\g\\
&=&\sum_{i=1}^{n-1}\add_x(\alpha(y_1),\cdots,\alpha(y_{i-1}),\add_x(x_1,\cdots,x_{n-2},y_i),\alpha(y_{i+1}),\cdots,\alpha(y_{n})),
\end{eqnarray*}
which implies that Conditions (iii) and (iv) in Definition \ref{derivation} hold. The proof is finished. \qed\vspace{3mm}

For any linear map $D:\wedge^{n-1}\g\lon\g$, denote by $\mathbb KD$ the $1$-dimensional vector space generated by $D$. On the direct sum $\g\oplus\mathbb KD$, define a totally skew-symmetric linear map $[\cdot,\cdots,\cdot]_D:\wedge^{n}(\g\oplus\mathbb KD)\lon\g\oplus\mathbb KD$ by
\begin{eqnarray*}
[x_1+k_1D,\cdots,x_n+k_nD]_D=[x_1,\cdots,x_n]_\g+\sum_{i=1}^n(-1)^{i-1}k_iD(x_1,\cdots,\hat{x}_i,\cdots,x_n).
\end{eqnarray*}
Define a linear map $\alpha_D:\g\oplus\mathbb KD\lon\g\oplus\mathbb KD$ by $\alpha_D(x+kD)=\alpha(x)+kD$, i.e.
\[
\alpha_D=\left(\begin{array}{cc}
\alpha & 0\\
0      & \Id
\end{array}\right).
\]

\begin{thm}
Let $(\g,[\cdot,\cdots,\cdot]_\g,\alpha)$ be an $n$-Hom-Lie algebra and $D:\wedge^{n-1}\g\lon\g$ a linear map. Then $(\g\oplus\mathbb KD,[\cdot,\cdots,\cdot]_D,\alpha_D)$ is an $n$-Hom-Lie algebra if and only if $D$ is generalized derivation on $\g$. We call $(\g\oplus\mathbb KD,[\cdot,\cdots,\cdot]_D,\alpha_D)$  the generalized derivation extension of $\g$ by the generalized derivation $D$.
\end{thm}

\pf For all $x_1,\cdots,x_n\in\g,k_1,\cdots,k_n\in\mathbb K$, we have
\begin{eqnarray*}
\alpha_D[x_1+k_1D,\cdots,x_n+k_nD]_D&=&\alpha_D([x_1,\cdots,x_n]_\g+\sum_{i=1}^n(-1)^{i-1}k_iD(x_1,\cdots,\hat{x}_i,\cdots,x_n))\\
                                    &=&\alpha([x_1,\cdots,x_n]_\g)+\sum_{i=1}^n(-1)^{i-1}k_i\alpha(D(x_1,\cdots,\hat{x}_i,\cdots,x_n)),
\end{eqnarray*}
and
\begin{eqnarray*}
&&[\alpha_D(x_1+k_1D),\cdots,\alpha_D(x_n+k_nD)]_D\\
&=&[\alpha(x_1)+k_1D,\cdots,\alpha(x_n)+k_nD]_D\\
                                                &=&[\alpha(x_1),\cdots,\alpha(x_n)]_\g+\sum_{i=1}^n(-1)^{i-1}k_iD(\alpha(x_1),\cdots,\hat{\alpha}(x_i),\cdots,\alpha(x_n)).
\end{eqnarray*}
Since $\alpha$ is an algebra morphism,   $\alpha_D$ is an algebra morphism if and only if
$$\alpha\circ D=D\circ\widetilde{\alpha}.$$
By the definition of the bracket $[\cdot,\cdots,\cdot]_D$, we deduce that $(\g\oplus\mathbb KD,[\cdot,\cdots,\cdot]_D,\alpha_D)$ satisfy the   Hom Fundamental identity if and only if
\begin{eqnarray}
&&\label{g-derivation1}\HF_{x_1,\cdots,x_{n-2},x_{n-1},D,y_1,\cdots,y_{n-1}}=0,\\
&&\label{g-derivation2}\HF_{D,x_1,\cdots,x_{n-2},y_1,\cdots,y_{n-1},y_n}=0,\\
&&\label{g-derivation3}\HF_{D,x_1,\cdots,x_{n-2},D,y_1,\cdots,y_{n-1}}=0.
\end{eqnarray}
By straightforward computation,  we have
\begin{eqnarray*}
 [\alpha_D(x_1),\cdots,\alpha_D(x_{n-1}),[D,y_1,\cdots,y_{n-1}]_D]_D
&=&[\alpha(x_1),\cdots,\alpha(x_{n-1}),D(y_1,\cdots,y_{n-1})]_D\\
&=&[\alpha(x_1),\cdots,\alpha(x_{n-1}),D(y_1,\cdots,y_{n-1})]_\g,
\end{eqnarray*}
and
\begin{eqnarray*}
&&\sum_{i=1}^{n-1}[\alpha_D(D),\alpha_D(y_1),\cdots,\alpha_D(y_{i-1}),
     [x_1,\cdots,x_{n-1},y_i]_D,\alpha_D(y_{i+1}),\cdots,\alpha_D(y_{n-1})]_D\\
&&+[[x_1,\cdots,x_{n-1},D]_D,\alpha_D(y_{1}),\cdots,\alpha_D(y_{n-1})]_D\\
&=&\sum_{i=1}^{n-1}[D,\alpha(y_1),\cdots,\alpha(y_{i-1}),
     [x_1,\cdots,x_{n-1},y_i]_\g,\alpha(y_{i+1}),\cdots,\alpha(y_{n-1})]_D\\
&&+[(-1)^{n-1}D(x_1,\cdots,x_{n-1}),\alpha(y_{1}),\cdots,\alpha(y_{n-1})]_D\\
&=&\sum_{i=1}^{n-1}D(\alpha(y_1),\cdots,\alpha(y_{i-1}),
     [x_1,\cdots,x_{n-1},y_i]_\g,\alpha(y_{i+1}),\cdots,\alpha(y_{n-1}))\\
&&+[\alpha(y_{1}),\cdots,\alpha(y_{n-1}),D(x_1,\cdots,x_{n-1})]_\g.
\end{eqnarray*}
Thus,   \eqref{g-derivation1} is equivalent to Condition (ii) in Definition \ref{derivation}.

Similarly, we have
\begin{eqnarray*}
[\alpha_D(D),\alpha_D(x_1),\cdots,\alpha_D(x_{n-2}),[y_1,\cdots,y_{n}]_D]_D
&=&[D,\alpha(x_1),\cdots,\alpha(x_{n-2}),[y_1,\cdots,y_{n}]_\g]_D\\
&=&D(\alpha(x_1),\cdots,\alpha(x_{n-2}),[y_1,\cdots,y_{n}]_\g),
\end{eqnarray*}
and
\begin{eqnarray*}
&&\sum_{i=1}^{n}[\alpha_D(y_1),\cdots,\alpha_D(y_{i-1}),
     [D,x_1,\cdots,x_{n-2},y_i]_D,\alpha_D(y_{i+1}),\cdots,\alpha_D(y_{n})]_D\\
&&=\sum_{i=1}^{n}[\alpha(y_1),\cdots,\alpha(y_{i-1}),
     D(x_1,\cdots,x_{n-2},y_i),\alpha(y_{i+1}),\cdots,\alpha(y_{n})]_D\\
&&=\sum_{i=1}^{n}[\alpha(y_1),\cdots,\alpha(y_{i-1}),
     D(x_1,\cdots,x_{n-2},y_i),\alpha(y_{i+1}),\cdots,\alpha(y_{n})]_\g.
\end{eqnarray*}
Thus,   \eqref{g-derivation2} is equivalent to Condition (iii) in Definition \ref{derivation}.

Finally, we have
\begin{eqnarray*}
&&[\alpha_D(D),\alpha_D(x_1),\cdots,\alpha_D(x_{n-2}),[D,y_1,\cdots,y_{n-1}]_D]_D\\
&&=[D,\alpha(x_1),\cdots,\alpha(x_{n-2}),D(y_1,\cdots,y_{n-1})]_D\\
&&=D(\alpha(x_1),\cdots,\alpha(x_{n-2}),D(y_1,\cdots,y_{n-1})),
\end{eqnarray*}
and
\begin{eqnarray*}
&&\sum_{i=1}^{n-1}[\alpha_D(D),\alpha_D(y_1),\cdots,\alpha_D(y_{i-1}),
     [D,x_1,\cdots,x_{n-2},y_i]_D,\alpha_D(y_{i+1}),\cdots,\alpha_D(y_{n-1})]_D\\
&&=\sum_{i=1}^{n-1}[D,\alpha(y_1),\cdots,\alpha(y_{i-1}),
     D(x_1,\cdots,x_{n-2},y_i),\alpha(y_{i+1}),\cdots,\alpha(y_{n-1})]_D\\
&&=\sum_{i=1}^{n-1}D(\alpha(y_1),\cdots,\alpha(y_{i-1}),
     D(x_1,\cdots,x_{n-2},y_i),\alpha(y_{i+1}),\cdots,\alpha(y_{n-1})).
\end{eqnarray*}
Thus,   \eqref{g-derivation3} is equivalent to Condition (iv) in Definition \ref{derivation}. Therefore, $(\g\oplus\mathbb KD,[\cdot,\cdots,\cdot]_D,\alpha_D)$ is an $n$-Hom-Lie algebra if and only if $D$ is a generalized derivation on $\g$. The proof is finished. \qed

\begin{pro}
Let $D^2$ and $D^1$ be two generalized derivations on an $n$-Hom-Lie algebra $(\g,[\cdot,\cdots,\cdot]_\g,\alpha)$. If there exists $x\in\g$ such that $\alpha(x)=x$ and $D^1=D^2+\add_x$, then the corresponding generalized derivation extensions $(\g\oplus\mathbb KD^2,[\cdot,\cdots,\cdot]_{D^2},\alpha_{D^2})$ and $(\g\oplus\mathbb KD^1,[\cdot,\cdots,\cdot]_{D^1},\alpha_{D^1})$ are isomorphic.
\end{pro}
\pf Define $\overline{x}:\K D^1\longrightarrow\g$ by
$$
\overline{x}(kD^1)=kx,\quad \forall k\in\K.
$$
Then $\left(\begin{array}{cc}\Id_\g&\overline{x}\\0&1\end{array}\right)$ is an isomorphism from the $n$-Hom-Lie algebra $(\g\oplus\mathbb KD^1,[\cdot,\cdots,\cdot]_{D^1},\alpha_{D^1})$ to $(\g\oplus\mathbb KD^2,[\cdot,\cdots,\cdot]_{D^2},\alpha_{D^2})$. We leave details to readers. \qed

\end{document}